\documentclass{amsart}
\usepackage[latin1]{inputenc}

\theoremstyle{plain}
\newtheorem{lemma}{Lemma}
\newtheorem{proposition}[lemma]{Proposition}
\newtheorem{theorem}[lemma]{Theorem}

\theoremstyle{definition}

\theoremstyle{remark}

\newcommand{\txt}[1]{\quad\mbox{#1}\quad}

\begin{document}

\title{Strong martingale type and uniform smoothness}
\author{Jörg Wenzel}
\date{\today}



\address{Jörg Wenzel, Department of Mathematics and Applied
  Mathematics, University of Pretoria, Pretoria 0002, South Africa}

\email{\url{wenzel@minet.uni-jena.de}}

\begin{abstract}
  We introduce stronger versions of the usual notions of martingale type
  $p\leq 2$ and cotype $q\geq 2$ of a Banach space $X$ and show that
  these concepts are equivalent to uniform $p$-smoothness and
  $q$-convexity, respectively. All these are metric concepts, so
  they depend on the particular norm in $X$.

  These concepts allow us to get some more insight into the fine line
  between $X$ being isomorphic to a uniformly $p$-smooth space or
  being uniformly $p$-smooth itself.

  Instead of looking at Banach spaces, we consider linear operators
  between Banach spaces right away. The situation of a Banach space
  $X$ can be rediscovered from this by considering the identity map of
  $X$.
\end{abstract}

\maketitle

\section{Introduction}
\label{sec:introduction}

In several recent papers Kato {\it et
  al.}~\cite{MR2002e:46015,MR2003e:46013} introduced the concept of
strong Rademacher type and cotype and related those concepts to the
uniform smoothness and convexity properties of the underlying Banach
space.

In particular for $1\leq p\leq 2$, a Banach space $X$ has \emph{strong
  Rademacher type $p$} if there is a constant $c$ such that
\[ \Big\| \sum_{k=0}^n x_kr_k\Big| L_p\Big\| \leq \Big(\|x_0\|^p +
c^p \sum_{k=1}^n \|x_k\|^p \Big)^{1/p}
\]
for all $n\in\mathbb{N}$ and elements $x_0,\dots,x_n\in X$. Here
$r_0,\dots,r_n$ denotes the sequence of Rademacher functions.

It turns out~\cite[Theorem 3]{MR2003e:46013} that $X$ is uniformly
$p$-smooth if and only if it has strong Rademacher type $p$, moreover
the constants involved coincide (see Section~\ref{sec:notation} for
the definition of $p$-smoothness and $q$-convexity).

A dual result holds for strong Rademacher cotype $q$ and uniform
$q$-convexity, where $q\geq 2$.

Surely, strong Rademacher type $p$ implies usual Rademacher type $p$
and the same goes for the respective cotype properties. Also, both
convexity and smoothness as well as the strong Rademacher type and
cotype properties depend on the particular norm in $X$.

It is well known~\cite{pietsch98:_orthon_banac,pis75}, that $X$ is
uniformly $p$-smooth for some equivalent norm, if and only if $X$ has
martingale type $p$. So the more natural notion connecting to uniform
$p$-smoothness would be something like \emph{strong martingale type
  $p$}.

In this paper, we define strong martingale type $p$ and show, that it
is equivalent (with the same constants) to $X$ being uniformly
$p$-smooth. Again, strong martingale type $p$ implies usual martingale
type $p$ but now we also have a partial converse. It follows, that one
can renorm a Banach space with martingale type $p$ so as to have
strong martingale type $p$ in the new norm. This is not so if one only
considers Rademacher type, since there is an example, due to
James~\cite{jam77} of a nonreflexive Banach space of Rademacher type
$2$. Since the space is non reflexive, it cannot be uniformly
$p$-smooth for any $p>1$ and consequently not have strong Rademacher
type $p$ for any $p>1$ in any equivalent norm.

All results have dual variants for strong martingale cotype $q$ and
uniform $q$-convexity.

Hardly any of the results in this paper are really new, they are just
old ones and generalizations of old ones, put into the right
perspective. Moreover, our direct proof for the renorming result in
the martingale type $p$ case (Theorem~\ref{thr:4}) seems to be new.

\section{Notation and concepts}
\label{sec:notation}

For $k=0,1,\dots$, the \emph{dyadic intervals}
\[ \Delta_k^{(i)}:=\Big[\frac {i}{2^k},\frac {i+1}{2^k}\Big)
\txt{where $i=0,\dots,2^k-1$,}
\]
generate the dyadic $\sigma$-algebra denoted by $\mathcal{F}_k$.

For a Banach space $X$, we consider \emph{dyadic martingales}
$(f_0,\dots,f_n)$, defined on $[0,1)$, taking values in $X$, and
adapted to the \emph{dyadic filtration}
$\mathcal{F}_0\subseteq\dots\subseteq\mathcal{F}_n$. We let
$f_{-1}\equiv0$ and denote by $d_k:=f_k-f_{k-1}$ the
\emph{differences} or increments of this martingale.

By
\[ \|f|L_p\| := \Big(\int_0^1\|f(t)\|^p\,dt\Big)^{1/p}
\]
we denote the $L_p$-norm of a function $f:[0,1)\to X$. Given another
function $f':[0,1)\to X'$, we write
\[ \langle f,f'\rangle := \int_0^1 \langle f(t),f'(t)\rangle\,dt,
\]
where $\langle x,x'\rangle$ is used to denote the duality between $X$
and $X'$.

For $q\geq 2$, an operator $T:X\to Y$ is \emph{uniformly $q$-convex}
if there is a constant $c$ such that
\begin{equation}
  \label{eq:1}
  \Big\| \frac{Tx_+-Tx_-}2\Big\| \leq c \bigg(
  \frac{\|x_+\|^q+\|x_-\|^q}2 -\Big\|
  \frac{x_++x_-}2\Big\|^q\bigg)^{1/q}
\end{equation}
for all $x_+,x_-\in X$. Equivalently, letting $x=(x_+-x_-)/2$ and
$x_0=(x_++x_-)/2$, we can rephrase the definition as
\begin{equation}
  \label{eq:12}
  \|Tx\| \leq c \Big( \frac{\|x+x_0\|^q + \|x-x_0\|^q}2
  -\|x_0\|^p\Big)^{1/q}
  \tag{1'}
\end{equation}
for all $x,x_0\in X$.

The operator has \emph{strong martingale cotype} $q$ if there is a
constant $c$ such that
\begin{equation}
  \label{eq:2}
  \Big(\|x\|^q + \frac1{c^q} \sum_{k=1}^n \|Td_k|L_q\|^q\Big)^{1/q}
  \leq
  \Big\|\sum_{k=0}^n d_k\Big|L_q\Big\|
\end{equation}
for all sequences of $X$-valued differences $d_0,d_1,\dots,d_n$ of
dyadic martingales such that $d_0\equiv x$.

Note that following Pisier~\cite{pis75} an operator has
\emph{martingale cotype $q$} if there is a constant $c$ such that
\begin{equation}\label{eq:7}
  \Big(\frac1{c^q} \sum_{k=0}^n \|Td_k|L_q\|^q\Big)^{1/q}
  \leq
  \Big\|\sum_{k=0}^n d_k\Big|L_q\Big\|
\end{equation}
for all sequences of $X$-valued differences $d_1,\dots,d_n$ of dyadic
martingales. So by $\|T\|\leq c$ strong martingale cotype $q$ implies
martingale cotype $q$. We will see, that the reverse is only true if
one passes to an equivalent norm in $X$.

We now repeat the definitions above in the dual case. For $1\leq p\leq
2$, an operator $T:X\to Y$ is \emph{uniformly $p$-smooth} if there is
a constant $c$ such that
\begin{equation}\label{eq:5}
  \bigg( \frac{\|y+Tx\|^p + \|y-Tx\|^p}2 -\|y\|^p\bigg)^{1/p} \leq c
  \|x\|
\end{equation}
for all $x\in X$ and $y\in Y$.

The operator has \emph{strong martingale type} $p$ if there is a
constant $c$ such that
\begin{equation}\label{eq:4}
  \Big\|y+\sum_{k=1}^n Td_k\Big|L_p\Big\| \leq \Big(\|y\|^p + c^p
  \sum_{k=1}^n \|d_k|L_p\|^p \Big)^{1/p}
\end{equation}
for all sequences of $X$-valued differences $d_1,\dots,d_n$ of
dyadic martingales and $y\in Y$.

Note that following Pisier~\cite{pis75} an operator has
\emph{martingale type $p$} if there is a constant $c$ such that
\begin{equation}\label{eq:8}
  \Big\|\sum_{k=0}^n Td_k\Big|L_p\Big\| \leq \Big(c^p
  \sum_{k=0}^n \|d_k|L_p\|^p \Big)^{1/p}
\end{equation}
for all sequences of $X$-valued differences $d_0,\dots,d_n$ of dyadic
martingales. So by letting $y=Tx$ and since $\|T\|\leq c$, strong
martingale type $p$ implies martingale type~$p$. We will see, that the
reverse is only true if one passes to an equivalent norm in~$Y$.

For $1\leq p <\infty$, we denote by $p'$ the \emph{dual index} of $p$,
given by $1/p+1/p'=1$.

\section{Results}
\label{sec:results}

\begin{theorem}\label{thr:1}
  An operator $T$ is uniformly $q$-convex with constant $c$ if and
  only if $T$ has strong martingale cotype $q$ with the same constant
  $c$.
\end{theorem}
\begin{proof}
  Assume first that $T$ has strong martingale cotype $q$. Define a
  martingale $d_0+d_1$ by letting $d_0\equiv x_0$ and $d_1:=x_1r_1$.
  Applying~(\ref{eq:2}) to this martingale yields
  \[ \Big( \|x_0\|^q + \frac1{c^q} \|Tx_1\|^q\Big)^{1/q}
  \leq \Big( \frac{\|x_0+x_1\|^q + \|x_0-x_1\|^q}2 \Big)^{1/q}.
  \]
  The substitution $x_+:=x_0+x_1$ and $x_-:=x_0-x_1$ then
  yields~(\ref{eq:1}).

  Conversely, if $d_0,\dots,d_n$ is any sequence of differences of a
  dyadic $X$-valued martingale such that $d_0\equiv x$, write
  $f_k=\sum_{h=0}^k d_h$. Then $f_{k-1}$ is constant on the intervals
  $\Delta_{k-1}^{(j)}$, while the (constant) value that $d_k$ takes on
  the left half $\Delta_k^{(2j)}$ of $\Delta_{k-1}^{(j)}$ is the
  negative of the (constant) value that $d_k$ takes on
  the right half $\Delta_k^{(2j+1)}$ of $\Delta_{k-1}^{(j)}$. Hence
  \[ \int\limits_{\Delta_{k-1}^{(j)}} \| f_{k-1}(t)+d_k(t)
  \|^q\, dt
  =
  \int\limits_{\Delta_{k-1}^{(j)}} \| f_{k-1}(t)  -
  d_k(t)\|^q\, dt
  \]
  and we have
  \begin{equation}
    \label{eq:3}
    \big\| f_k \big|L_q \big\| = \big\| f_{k-1}
    + d_k\big|L_q \big\| = \big\| f_{k-1}
    - d_k\big|L_q \big\|.
  \end{equation}
  From~(\ref{eq:12}) we get
  \[ \|Td_k(t)\|^q \leq c^q \Big( \frac{\|f_{k-1}(t) +
    d_k(t)\|^q + \|f_{k-1}(t) - d_k(t)\|^q}2 -
  \big\|f_{k-1}(t)\big\|^q \Big)
  \]
  which, when integrated over $t$ and using~(\ref{eq:3}) gives
  \[ \|Td_k|L_q\|^q \leq c^q \big( \|f_k|L_q\|^q -
  \|f_{k-1}|L_q\|^q \big).
  \]
  Take the sum over $k=1,\dots,n$ of these inequalities to
  get~(\ref{eq:2}).
\end{proof}

\begin{theorem}\label{thr:2}
  An operator $T$ is uniformly $p$-smooth with constant $c$ if and
  only if $T$ has strong martingale type $p$ with the same constant
  $c$.
\end{theorem}
\begin{proof}
  Assume first that $T$ has strong martingale type $p$. Define a
  martingale difference by letting $d_1\equiv xr_1$.
  Applying~(\ref{eq:4}) to this martingale yields
  \[ \Big( \frac{\|y+Tx\|^p+\|y-Tx\|^p}2 \Big)^{1/p} \leq \big(
  \|y\|^p + c^p \|x\|^p \big)^{1/p},
  \]
  which immediately yields~(\ref{eq:5}).

  Conversely, if $d_1,\dots,d_n$ is any sequence of differences of a
  dyadic $X$-valued martingale, write $f_k=\sum_{h=1}^n d_h$. Then
  $y+Tf_{k-1}$ is constant on the intervals $\Delta_{k-1}^{(j)}$,
  while the (constant) value that $Td_k$ takes on the left half
  $\Delta_k^{(2j)}$ of $\Delta_{k-1}^{(j)}$ is the negative of the
  (constant) value that $Td_k$ takes on the right half
  $\Delta_k^{(2j+1)}$ of $\Delta_{k-1}^{(j)}$. Hence
  \[ \int\limits_{\Delta_{k-1}^{(j)}} \|y+Tf_{k-1}(t)+d_k(t)\|^p\, dt
  =  \int\limits_{\Delta_{k-1}^{(j)}} \|y+Tf_{k-1}(t)-d_k(t)\|^p\, dt
  \]
  and we have
  \begin{equation}
    \|y+Tf_k|L_p\| = \|y+Tf_{k-1} + Td_k|L_p\| = \|y+ Tf_{k-1} -
    Td_k|L_p\|.
    \label{eq:6}
  \end{equation}
  From~(\ref{eq:5}) we get
  \begin{multline*}
    \frac{\|y+Tf_{k-1}(t) + Td_k(t)\|^p + \|y+Tf_{k-1}(t) -
      Td_k(t)\|^p} 2\\
    \leq \|y+Tf_{k-1}(t)\|^p+ c^p \|d_k(t)\|^p
  \end{multline*}
  which when integrated over $t$ and using~(\ref{eq:6}) gives
  \[ \|y+Tf_k|L_p\|^p \leq \|y+Tf_{k-1}|L_p\|^p + c^p\|d_k|L_p\|^p.
  \]
  Take the sum over $k=1,\dots,n$ of these inequalities to
  get~(\ref{eq:4}).
\end{proof}

\begin{proposition}\label{prop:2}
  An operator $T$ has strong martingale type $p$ if and only if its
  dual $T'$ has strong martingale cotype $p'$. And an operator $T$ has
  strong martingale cotype $q$ if and only if its dual $T'$ has strong
  martingale type $q'$.
\end{proposition}
\begin{proof}
  Assume that $T'$ is of strong martingale cotype $p'$.

  Given a sequence $d_1,\dots,d_n$ of dyadic martinagle differences
  and $y\in Y$, we find an $\mathcal{F}_n$-measurable function $g_n$
  with values in $Y'$ such that
  \[ \Big\| y + \sum_{k=1}^n Td_k\Big| L_p\Big\| = \Big\langle
  y+ \sum_{k=1}^n Td_k,g_n\Big\rangle
  \txt{and} \|g_n|L_{p'}\|\leq 1+\epsilon.
  \]
  Writing
  \[ e_k := \mathbb{E}(g_n|\mathcal{F}_k) -
  \mathbb{E}(g_n|\mathcal{F}_{k-1}),
  \]
  we obtain another sequence of dyadic martingale differences and we
  have
  \[ \Big\langle y+
  \sum_{k=1}^n Td_k,g_n\Big\rangle = \langle y,e_0\rangle +
  \sum_{k=1}^n \langle Td_k,e_k\rangle.
  \]
  Consequently
  \begin{multline*}
    \Big\| y+ \sum_{k=1}^n Td_k\Big| L_p\Big\| \leq \\
    \Big(\|y\|^p
    + c^p \sum_{k=1}^n \|d_k|L_p\|^p \Big)^{\!\!1/p} \Big(
    \|e_0\|^{p'} + \frac1{c^{p'}} \sum_{k=1}^n
    \|T'e_k|L_{p'}\|^{p'} \Big)^{\!\!1/p'}.
  \end{multline*}
  Now applying the strong martingale cotype $p'$ property of $T'$ we
  get
  \[ \|e_0\|^{p'} + \frac1{c^{p'}} \sum_{k=1}^n
  \|T'e_k|L_{p'}\|^{p'}
  \leq \Big\|
  \sum_{k=0}^n e_k\Big| L_{p'}\Big\|^{p'} \leq
  (1+\epsilon)^{p'},
  \]
  which proves that
  \[ \Big\| y+ \sum_{k=1}^n Td_k\Big| L_p\Big\| \leq
  (1+\epsilon) \Big( \|y\|^p + c^p \sum_{k=1}^n \|d_k|L_p\|^p
  \Big)^{1/p}.
  \]
  Letting $\epsilon\to 0$ proves that $T$ has strong martingale type
  $p$ with the same constant.

  To see the other implication, given a sequence of $Y'$-valued dyadic
  martingale differences $d_0,\dots,d_n$ such that $d_0\equiv y'$, let
  \[ \lambda := \frac{c\|y'\|^{p'-1}}{\displaystyle
    \bigg( \Big\| \sum_{k=0}^n
    d_k \Big| L_{p'}\Big\|^{p'} - \|y'\|^{p'}\bigg)^{1/p}}.
  \]
  We find $X'$-valued dyadic martingale differences $e_1,\dots,e_n$
  such that
  \[ \Big(\sum_{k=1}^n \|T'd_k|L_{p'}\|^{p'}\Big)^{1/p'} =
  \sum_{k=1}^n
  \langle T'd_k,e_k\rangle \txt{and} \Big(\sum_{k=1}^n
  \|e_k|L_p\|^p\Big)^{1/p} \leq 1+\epsilon.
  \]
  Moreover we find $y\in Y$ such that $\|y'\|=\langle y',y \rangle$
  and $\|y\|\leq 1+\epsilon$.

  We can then write
  \begin{align*}
    \Big(\sum_{k=1}^n \|T'd_k|L_{p'}\|^{p'}\Big)^{1/p'}
    &=
    \Big\langle y' + \sum_{k=1}^n d_k,\lambda y+\sum_{k=1}^n
    Te_k\Big\rangle - \lambda\langle y',y\rangle \\
    &\leq
    \Big\| \sum_{k=0}^n d_k \Big| L_{p'}\Big\|\cdot \Big\| \lambda
    y + \sum_{k=1}^n Te_k\Big|L_p\Big\| - \lambda\|y'\|.
  \end{align*}
  Using the strong martingale type $p$ property of $T$, we get
  \[ \Big\| \lambda y + \sum_{k=1}^n Te_k\Big|L_p\Big\|
  \leq \Big( \lambda^p \|y\|^p + c^p \sum_{k=1}^n \|e_k|L_p\|^p
  \Big)^{1/p} \leq (1+\epsilon) \big( c^p +\lambda^p \big)^{1/p}.
  \]
  Our choice of $\lambda$ now ensures that
  \begin{equation*}
    \Big\| \sum_{k=0}^n d_k \Big| L_{p'}\Big\| ( c^p
    + \lambda^p )^{1/p}
    = c\bigg( \Big\|\sum_{k=0}^n d_k\Big| L_{p'}\Big\|^{p'}
    - \|y'\|^{p'}\bigg)^{1/p'} + \lambda\|y'\|.
  \end{equation*}
  This proves that
  \begin{align*}
    \Big( \sum_{k=1}^n \|T'd_k|L_{p'}\|^{p'}\Big)^{1/p'}
    &\leq
    \Big\|\sum_{k=0}^n d_k\Big| L_{p'}\Big\|
    (c^p+\lambda^p)^{1/p}(1+\epsilon) - \lambda\|y'\| \\
    &=
    (1+\epsilon)c \bigg( \Big\| \sum_{k=0}^n d_k\Big| L_{p'}\Big\|^{p'}
    - \|y'\|^{p'} \bigg)^{1/p'} + \epsilon\lambda\|y'\|.
  \end{align*}
  Letting $\epsilon\to0$ proves that $T'$ has strong martingale cotype
  with the same constant.

  The second part follows analogously.
\end{proof}

\begin{proposition}\label{prop:1}
  An operator $T$ is uniformly $p$-smooth if and only if its dual $T'$
  is uniformly $p'$-convex. An operator $T$ is uniformly $q$-convex if
  and only if its dual $T'$ is uniformly $q'$-smooth.
\end{proposition}
\begin{proof}
  This was well known for Banach spaces, see
  Beauzamy~\cite[pp.~311--312]{bea} and has been proved for linear
  operators in Pietsch\slash
  Wenzel~\cite[7.9.6]{pietsch98:_orthon_banac}.
\end{proof}

We have now an alternative proof of Theorem~\ref{thr:2} using duality.

\begin{proof}[Proof of Theorem~\ref{thr:2} (alternative version)]
  If $T$ is uniformly $p$-smooth, then its dual $T'$ is uniformly
  $p'$-convex, which by Theorem~\ref{thr:1} happens if and only if
  $T'$ has strong martingale cotype $p'$ and by
  Proposition~\ref{prop:2} this is equivalent to $T$ having strong
  martingale type~$p$.
\end{proof}

As our last project, we want to discuss the relation of martingale
type\slash cotype and strong martingale type\slash cotype. It is due
to Pisier~\cite{pis75} that whenever $T$ has martingale cotype $q$,
then there exists an equivalent norm on $X$, such that $T$ becomes
uniformly $q$-convex when considered as an operator from $X$ equipped
with this new norm into $Y$. From this it then follows that $T$ has
strong martingale cotype $q$. By duality the same works in the type
case However, we also want to give a direct argument similar to the
one used for cotype.

Since we think that our approach makes more clear, why the respective
proofs work and since we haven't seen the direct proof for the type
case in print yet, we want to include a proof for both cases here.

We first provide a technical lemma, that allows us to define
equivalent norms.
\begin{lemma}\label{lem:1}
  Let $|||\,\cdot\,|||:[X,\|\,\cdot\,\|]\to [0,\infty)$ be a
  continous, positively homogeneous functional. Assume that for
  $|||x_\pm|||\leq 1$ it follows that $|||x_++x_-|||\leq 2$. Then
  \[ |||x_++x_-|||\leq |||x_+|||+|||x_-|||
  \]
  for all $x_\pm\in X$.
\end{lemma}
\begin{proof}
  Assume first that $|||x_\pm|||\leq 1$. By assumption
  \begin{equation}
    |||\lambda x_++(1-\lambda)x_-|||\leq 1
    \label{eq:11}
  \end{equation}
  for $\lambda=1/2$. It follows by induction over $n$, that the same
  holds for all $\lambda=k/2^n$, where $n=0,1,2,\dots$ and $0\leq
  k\leq 2^n$.

  To see this, write
  \[ x:=\frac k{2^{n-1}} x_+ + \Big(1-\frac k{2^{n-1}}\Big) x_-
  \]
  and note that $|||x|||\leq 1$ by the induction hypothesis, and
  \[ \lambda x_+ + (1-\lambda)x_- = \frac12x+\frac12x_-.
  \]

  It now follows by continuity that~(\ref{eq:11}) holds in fact for
  all $\lambda\in[0,1]$.

  Finally, for arbitrary $x_\pm$ let
  \[ \lambda=\frac{|||x_+|||}{|||x_+|||+|||x_-|||}.
  \]
  Then from~(\ref{eq:11}) we get
  \[ \Big|\Big|\Big| \lambda \frac{x_+}{|||x_+|||} + (1-\lambda)
  \frac{x_-}{|||x_-|||} \Big|\Big|\Big| \leq  1
  \]
  which in turn implies the triangle inequality.
\end{proof}

\begin{theorem}\label{thr:3}
  If an operator $T:X\to Y$ has martingale cotype $q$ then there
  exists an equivalent norm on $X$, such that $T$ considered as an
  operator from $X$ equipped with the new norm into $Y$ is uniformly
  $q$-convex.
\end{theorem}
\begin{proof}
  Assume that~(\ref{eq:7}) holds. Let
  \[ \{x\} := \inf \bigg( \Big\|x+\sum_{k=1}^n d_k\Big|L_q\Big\|^q -
  \frac1{c^q} \sum_{k=1}^n \|Td_k|L_q\|^q\bigg)^{1/q},
  \]
  where the infimum is taken over all sequences of $X$-valued
  differences $d_1,\dots,d_n$ of martingales adapted to the dyadic
  filtration $\mathcal{F}_1\subseteq\dots\subseteq\mathcal{F}_n$.

  Letting $d_1\equiv\dots\equiv d_n\equiv 0$ yields
  \[ \{x\}\leq \|x\|.
  \]
  Conversely, it follows from~(\ref{eq:7}) that
  \[ \{x\}\geq \frac1c \|Tx\|.
  \]

  Trivially with this expression (which need not yet be a norm) we
  have
  \[ \Big( \{x\}^q + \frac1{c^q} \sum_{k=1}^n \|Td_k|L_q\|^q
  \Big)^{1/q} \leq \Big\|x+\sum_{k=1}^n d_k\Big| L_q\Big\|
  \]
  which unfortunately is not yet strong martingale cotype $q$, since
  we would have to use $\{\,\cdot\,\}$ also on the right hand side.

  Therefore we next show that this expression satisfies~(\ref{eq:1}).
  Given $x_+$ and $x_-$ choose differences $d_1^\pm,\dots,d_n^\pm$ of
  dyadic martingales such that
  \[ \Big\| x_\pm+ \sum_{k=1}^n d_k^\pm \Big|L_q\Big\|^q - \frac1{c^q}
  \sum_{k=1}^n\|Td_k^\pm|L_q\|^q \leq \{x_{\pm}\}^q + \epsilon.
  \]
  Glueing together two differences as
  \[ d_{k+1}(t) :=
  \begin{cases}
    d_k^+(2t) & \mbox{if $0\leq t<1/2$,} \\
    d_k^-(2t-1) & \mbox{if $1/2\leq t<1$,}
  \end{cases}
  \]
  we get a new sequence of $X$-valued differences $d_2,\dots,d_{n+1}$
  of a dyadic martingale, which however is now adapted to
  $\mathcal{F}_2\subseteq\dots\subseteq\mathcal{F}_{n+1}$, so letting
  \[ d_1(t):=
  \begin{cases}
    \frac{x_+-x_-}2 & \mbox{if $0\leq t<1/2$,} \\
    \frac{x_--x_+}2 & \mbox{if $1/2\leq t<1$,}
  \end{cases}
  \]
  yields a sequence of differences of a dyadic martingale adapted to
  $\mathcal{F}_1\subseteq\dots\subseteq\mathcal{F}_{n+1}$.

  This is, by the way, the point, where we cannot just use Rademacher
  functions, since the function equal to $x_k^+r_k(2t)$ for $0 \leq
  t<1/2$ and to $x_k^-r_k(2t-1)$ for $1/2\leq t<1$ will no longer be a
  multiple of a Rademacher function, but at best a martingale
  difference as soon as $x_k^+\neq x_k^-$.

  Continuing with our considerations, we now have consequently
  \[ \{x\}^q \leq \Big\| x+\sum_{k=1}^{n+1} d_k\Big|L_q\Big\|^q -
  \frac1{c^q} \sum_{k=1}^{n+1} \|Td_k|L_q\|^q
  \]
  for all $x\in X$, in particular for $(x_++x_-)/2$.

  It is now clear that
  \[ \|Td_1|L_q\|^q = \Big\|\frac{Tx_+-Tx_-}2\Big\|^q, \quad
  \|Td_{k+1}|L_q\|^q = \frac{\|Td_k^+|L_q\|^q + \|Td_k^-|L_q\|^q}2,
  \]
  and
  \[ \Big\| \frac{x_++x_-}2 + \sum_{k=1}^{n+1}
  d_k\Big|L_q\Big\|^q = \frac12 \Big\| x_+ + \sum_{k=1}^n
  d_k^+\Big|L_q\Big\|^q + \frac12 \Big\|x_- + \sum_{k=1}^n
  d_k^-\Big|L_q\Big\|^q.
  \]
  Therefore it follows from the definition of $d_k^\pm$ that
  \[ \Big\{ \frac{x_++x_-}2\Big\}^q \leq
  \frac{\{x_+\}^q+\{x_-\}^q}2 +\epsilon -
  \frac1{c^q}\Big\|\frac{Tx_+-Tx_-}2\Big\|^q.
  \]
  Letting $\epsilon\to0$ shows
  \begin{equation}
    \label{eq:9}
    \Big\{ \frac{x_++x_-}2\Big\}^q \leq
    \frac{\{x_+\}^q+\{x_-\}^q}2 -
    \frac1{c^q}\Big\|\frac{Tx_+-Tx_-}2\Big\|^q.
  \end{equation}
  However, as was mentioned earlier, the expression $\{\,\cdot\,\}$
  need not be a norm on $X$. But it is positively homogeneous and
  continous so from~(\ref{eq:9}) it follows that Lemma~\ref{lem:1}
  applies and it also satisfies the triangle inequality. To get an
  \emph{equivalent} norm, we define
  \[ |||x||| := \big(\|x\|^q + \{x\}^q\big)^{1/q}.
  \]
  Adding~(\ref{eq:9}) to
  \[ \Big\| \frac{x_++x_-}2\Big\|^q \leq
  \frac{\|x_+\|^q+\|x_-\|^q}2
  \]
  we finally get the uniform $q$-convexity for $T$ considered as an
  operator from $X$ equipped with $|||\,\cdot\,|||$ to $Y$.
\end{proof}

\begin{theorem}\label{thr:4}
  If an operator $T:X\to Y$ has martingale type $p$ then there
  exists an equivalent norm on $Y$, such that $T$ considered as an
  operator from $X$ into $Y$ equipped with this new norm is uniformly
  $p$-smooth.
\end{theorem}
\begin{proof}
  Assume that~(\ref{eq:8}) holds. For $x\in X$ let
  \[ \{x\} := \sup \bigg( \Big\|Tx+\sum_{k=1}^n Td_k \Big|L_p\Big\|^p -
  c^p \sum_{k=1}^n\|d_k|L_p\|^p \bigg)^{1/p},
  \]
  where the supremum is taken over all sequences of $X$-valued
  differences $d_1,\dots,d_n$ of martingales adapted to the dyadic
  filtration $\mathcal{F}_1\subseteq\dots\subseteq\mathcal{F}_n$.

  Letting $d_1\equiv\dots\equiv d_n\equiv 0$ yields
  \[ \{x\} \geq \|Tx\|.
  \]
  Conversely, it follows from~(\ref{eq:8}) that
  \[ \{x\} \leq c \|x\|.
  \]
  We next show, that the expression $\{\,\cdot\,\}$
  satisfies
  \[ \Big\{ \frac{x_++x_-}2\Big\} ^p \geq \frac{\{x_+\}^p +
    \{x_-\}^p}2 - c^p \Big\|\frac{x_+-x_-}2\Big\|^p.
  \]
  Given $x_+$ and $x_-$ choose differences $d_1^\pm,\dots,d_n^\pm$ of
  dyadic martingales such that
  \[ \Big\| Tx_\pm+ \sum_{k=1}^n Td_k^\pm \Big|L_p\Big\|^p - c^p
  \sum_{k=1}^n \|d_k^\pm|L_p\|^p \geq \{x_\pm\} -\epsilon.
  \]
  Glueing together two differences as
  \[ d_{k+1}(t) :=
  \begin{cases}
    d_k^+(2t) & \mbox{if $0\leq t<1/2$,} \\
    d_k^-(2t-1) & \mbox{if $1/2\leq t<1$,}
  \end{cases}
  \]
  we get a new sequence of $X$-valued differences $d_2,\dots,d_{n+1}$
  of a dyadic martingale, which however is now adapted to
  $\mathcal{F}_2\subseteq\dots\subseteq\mathcal{F}_{n+1}$, so letting
  \[ d_1(t):=
  \begin{cases}
    \frac{x_+-x_-}2 & \mbox{if $0\leq t<1/2$,} \\
    \frac{x_--x_+}2 & \mbox{if $1/2\leq t<1$,}
  \end{cases}
  \]
  yields a sequence of differences of a dyadic martingale adapted to
  $\mathcal{F}_1\subseteq\dots\subseteq\mathcal{F}_{n+1}$.

  Consequently
  \[ \{x\}^p \geq \Big\| Tx+\sum_{k=1}^{n+1} Td_k\Big|L_p\Big\|^p -
  c^p \sum_{k=1}^{n+1} \|d_k|L_p\|^p
  \]
  for all $x\in X$, in particular for $(x_++x_-)/2$. It is now clear
  that
  \[ \|d_1|L_p\|^p = \Big\|\frac{x_+-x_-}2\Big\|^p, \quad
  \|d_k|L_p\|^p = \frac{\|d_k^+|L_p\|^p + \|d_k^-|L_p\|^p}2,
  \]
  and
  \[ \Big\| \frac{Tx_++Tx_-}2 + \sum_{k=1}^{n+1}
  Td_k\Big|L_p\Big\|^p = \frac12 \Big\| Tx_+ + \sum_{k=1}^n
  Td_k^+\Big|L_p\Big\|^p + \frac12 \Big\|Tx_- + \sum_{k=1}^n
  Td_k^-\Big|L_p\Big\|^p.
  \]
  Therefore it follows from the definition of $d_k^\pm$ that
  \[ \Big\{ \frac{x_++x_-}2\Big\}^p \geq \frac{\{x_+\}^p + \{x_-\}^p}2
  -\epsilon - c^p \Big\|\frac{x_+-x_-}2\Big\|^p.
  \]
  Letting $\epsilon\to0$ shows
  \begin{equation}
    \Big\{ \frac{x_++x_-}2\Big\}^p \geq \frac{\{x_+\}^p + \{x_-\}^p}2
    - c^p \Big\|\frac{x_+-x_-}2\Big\|^p.
    \label{eq:10}
  \end{equation}
  Finally defining
  \[ |||y||| := \inf \Big( \frac1{2^n} \sum_{k=1}^{2^n} \big(
  \|y_k-Tx_k\|^p + \{x_k\}^p \big) \Big)^{1/p}
  \]
  where the infimum is taken over all decompositions of $y$ as
  $y=\sum_{k=1}^{2^n} y_k/2^n$ and all elements $x_k\in X$, we we get
  an equivalent norm on $Y$ for which $T$ becomes uniformly
  $p$-smooth. Indeed clearly $|||\,\cdot\,|||$ is positively
  homogeneous and choosing $n=0$ and $x_0=0$ we get
  \[ |||y|||\leq \|y\|.
  \]
  On the other hand
  \begin{align*}
    \Big( \frac1{2^n} \sum_{k=1}^{2^n} \big( \|y_k-Tx_k\|^p +
    \{x\}^p \big) \Big)^{1/p}
    &\geq \Big( \frac1{2^n} \sum_{k=1}^{2^n}
    \big( \|y_k-Tx_k\|^p + \|Tx_k\|^p \big) \Big)^{1/p} \\
    &\geq \Big(
    \frac1{2^n} \sum_{k=1}^{2^n} 2^{1-p} \big( \|y_k-Tx_k\| + \|Tx_k\|
    \big)^p \Big)^{1/p} \\
    & \geq 2^{1/p-1} \Big( \frac1{2^n}
    \sum_{k=1}^{2^n} \|y_k\|^p \Big)^{1/p} \geq 2^{1/p-1} \|y\|
  \end{align*}
  whenever $y=\sum_{k=1}^{2^n}/2^n$, so
  \[ |||y|||\geq 2^{1/p-1} \|y\|.
  \]
  Next, given $y_\pm$ we choose $y_k^\pm$ and $x_k^\pm\in X$ such
  that
  \[ \frac 1{2^n} \sum_{k=1}^{2^n} \big(
  \|y_k^\pm-Tx_k^\pm\|^p + \{x_k^\pm\}^p \big) \leq |||y_\pm|||^p +
  \epsilon
  \txt{and}
  y_\pm=\frac1{2^n}\sum_{k=1}^{2^n} y_k^\pm.
  \]
  Note that we can assume that for both $y_+$ and $y_-$ the same $n$
  works, since if $y=\frac1{2^n} \sum_{k=1}^{2^n} y_k$ then also
  \[ y= \frac1{2^{n+1}} \Big( \sum_{k=1}^{2^n} y_k +
  \sum_{k=1}^{2^n} y_k \Big)
  \]
  and moreover, for $x_1,\dots,x_{2^n}\in X$ we have
  \begin{multline*}
    \bigg( \frac1{2^{n+1}} \Big( \sum_{k=1}^{2^n} \big( \|y_k-Tx_k\|^p
    + \{x_k\}^p \big) + \sum_{k=1}^{2^n} \big( \|y_k-Tx_k\|^p +
    \{x_k\}^p \big) \Big) \bigg)^{1/p} =\\
    \Big( \frac1{2^n}
    \sum_{k=1}^{2^n} \big( \|y_k-Tx_k\|^p + \{x_k\}^p\big) \Big)^{1/p}.
  \end{multline*}
  It now follows that
  \[ y_++y_- = \frac1{2^{n+1}} \Big( \sum_{k=1}^{2^n} (2y_k^+) +
  \sum_{k=1}^{2^n} (2y_k^-)\Big)
  \]
  hence
  \begin{align*}
    \lefteqn{ |||y_++y_-|||^p } \\
    & \leq \frac1{2^{n+1}} \Big( \sum_{k=1}^{2^n}
    \big( \|2y_k^+-2Tx_k^+\|^p + \{2x_k^+\}^p\big) +
    \sum_{k=1}^{2^n} \big(
    \|2y_k^--2Tx_k^-\|^p + \{2x_k^-\}^p \big) \Big)\\
    &= 2^{p-1} \Big( \frac1{2^n} \sum_{k=1}^{2^n} \big(
    \|y_k^+-Tx_k^+\|^p +\{x_k^+\}^p\big) + \frac1{2^n}
    \sum_{k=1}^{2^n} \big( \|y_k^--Tx_k^-\|^p + \{x_k^-\}^p\big) \Big)
    \\
    &\leq 2^{p-1} \big(|||y_+|||^p+|||y_-|||^p + 2\epsilon),
  \end{align*}
  which proves that
  \[ \Big|\Big|\Big| \frac{y_++y_-}2 \Big|\Big|\Big| \leq \Big(
  \frac{|||y_+|||^p+|||y_-|||^p}2 \Big)^{1/p}.
  \]
  It now follows from Lemma~\ref{lem:1} that $|||\,\cdot\,|||$ is
  actually a norm.

  To see the $p$-smoothness property of $T$ for this norm, given $y\in
  Y$, choose $y_k$ and $x_k\in X$ such that
  \[ \frac1{2^n}\sum_{k=1}^{2^n} \big( \|y_k-Tx_k\|^p +
  \{x_k\}^p\big)  \leq |||y|||^p + \epsilon
  \txt{and}
  y= \frac1{2^n}\sum_{k=1}^{2^n} y_k.
  \]
  Then for $x\in X$
  \begin{align*}
    |||y\pm Tx||| ^p
    &\leq \frac1{2^n} \sum_{k=1}^{2^n} \big( \|y_k\pm Tx
    - T(x_k\pm x)\|^p + \{x_k\pm x\}^p\big)\\
    &= \frac1{2^n} \sum_{k=1}^{2^n} \big( \|y_k
    - Tx_k\|^p + \{x_k\pm x\}^p\big).
  \end{align*}
  But by~(\ref{eq:10}) we have for each $k$
  \[ \frac{\{x_k+x\}^p + \{x_k-x\}^p}2 \leq \{x_k\}^p + c^p\|x\|^p
  \]
  so that
  \begin{align*}
    \frac{|||y + Tx||| ^p+|||y - Tx||| ^p}2
    &\leq \frac1{2^n}
    \sum_{k=1}^{2^n} \big( \|y_k-Tx_k\|^p +\{x_k\}^p + c^p \|x\|^p\big)
    \\
    &\leq |||y|||^p+\epsilon + c^p \|x\|^p
  \end{align*}
  which proves the uniform $p$-smoothness of $T$.
\end{proof}


\end{document}